\documentclass[11pt]{article}
\usepackage{amssymb}
\usepackage{mathptmx}

\input{tcilatex}

\def\widehat#1{\hat{#1}}

\begin{document}

\author{R\ L\ Hudson\thanks{%
School of Computing and Mathematics, Nottingham Trent University, Burton
Street, Nottingham NG1 4BU, UK. \texttt{robin.hudson@ntu.ac.uk}}
and S Pulmannov\'{a}\thanks{%
Mathematical Institute, Slovak Academy of Sciences, Stefanikova 49, 81473
Bratislava, Slovakia. \texttt{pulmann@mat.savba.sk}}}
\title{Double product integrals and Enriquez quantisation of Lie bialgebras II:
The quantum Yang-Baxter equation\thanks{%
AMS Classification 53D55, 17B62; keywords: quantisation, Lie bialgebras,
quantum Yang-Baxter equation}.}
\date{Received 30 September 2003, revised 21 April 2005}
\maketitle

\begin{abstract}
For a Lie algebra with Lie bracket got by taking commutators in a nonunital
associative algebra $\mathcal{L}$, let $\mathcal{T(L)}$ be the vector space
of tensors over $\mathcal{L}$ equipped with the It\^{o} Hopf algebra
structure derived from the associative multiplication in $\mathcal{L}$ $.$
It is shown that a necessary and sufficient condition that the double
product integral $\stackrel{\rightarrow \leftarrow }{\prod }(1+dr[h])$
satisfy the quantum Yang-Baxter equation over $\mathcal{T(L)}$ is that $%
1+r[h]$ satisfy the same equation over the unital associative algebra $%
\mathcal{L}^{\prime }$ got by adjoining a unit element to $\mathcal{L}$. In
particular the first-order coefficient $r_1$ of $r[h]$ satisfies the
classical Yang-Baxter equation. Using the fact that the multiplicative
inverse of $\stackrel{\rightarrow \leftarrow }{\prod }(1+dr[h])$ is $%
\stackrel{\leftarrow \rightarrow }{\prod }(1+dr^{\prime }[h])$ where $%
1+r^{\prime }[h]$ is the inverse of $1+r[h]$ in $\mathcal{L}^{\prime }[[h]]$
we construct a quantisation of an arbitrary quasitriangular Lie bialgebra
structure on $\mathcal{L}$ in the unital associative subalgebra of $\mathcal{%
T(L)}[[h]]$ consisting of formal power series whose zero order coefficient
lies in the space $\mathcal{S(L)}$ of symmetric tensors. The deformation
coproduct acts on $\mathcal{S(L)}$ by conjugating the undeformed coproduct
by $\stackrel{\rightarrow \leftarrow }{\prod }(1+dr[h])$ and the coboundary
structure $r$ of $\mathcal{L}$ is given by $r=r_1-\tau _{(2,1)}r_1$ where $%
\tau _{(2,1)}$ is the flip.
\end{abstract}

\section{\protect\medskip Introduction.}

Let $\mathcal{L}$ be a finite dimensional Lie algebra in which the Lie
bracket is got by taking commutators in a not necessarily unital associative
algebra $\mathcal{L}$. In a previous article \cite{HuPu} (see also \cite
{Enri1}) we showed that, in the Hopf algebra $\mathcal{T(L)}[[h]]$ of formal
power series with coefficients in the space $\mathcal{T(L)}$ of all tensors
over $\mathcal{L}$ equipped with the noncommutative It\^{o} extension of the
shuffle product, nonzero elements $R[h]$ of $\left( \mathcal{T(L)}\otimes
\mathcal{T(L)}\right) [[h]]$ satisfying the quasitriangularity relations
\begin{equation}
\left( \Delta \otimes \func{id}\,_{\mathcal{T(L)}}\right)
R[h]=R^{1,3}[h]R^{2,3}[h],\text{ }\left( \func{id}\,_{\mathcal{T(L)}}\otimes
\Delta \right) R[h]=R^{1,3}[h]R^{1,2}[h]  \label{triangular}
\end{equation}
\bigskip in the convolution algebra $\left( \mathcal{T(L)}\otimes \mathcal{%
T(L)}\otimes \mathcal{T(L)}\right) [[h]],$ where $\Delta $ is the coproduct,
can be characterised as directed double product integrals of form $\stackrel{%
\rightarrow \leftarrow }{\prod }(1+dr[h]),$ where $r[h]$ is an element of
the space $h\left( \mathcal{L}\otimes \mathcal{L}\right) [[h]]$ of formal
power series with coefficients in $\mathcal{L}\otimes \mathcal{L}$ and
vanishing constant term.

In the present article, we shall first establish the following necessary and
sufficient condition on $r[h]$ for such a product integral $R[h]=\stackrel{%
\rightarrow \leftarrow }{\prod }(1+dr[h])$ to satisfy the quantum
Yang-Baxter equation
\begin{equation}
R^{1,2}[h]R^{1,3}[h]R^{2,3}[h]=R^{2,3}[h]R^{1,3}[h]R^{1,2}[h]  \label{qybe}
\end{equation}
in the algebra of formal power series $\left( \mathcal{T(L)}\otimes \mathcal{%
T(L)}\otimes \mathcal{T(L)}\right) [[h]]$.

\begin{theorem}
A necessary and sufficient condition that the product integral $\stackrel{%
\rightarrow \leftarrow }{\prod }(1+dr[h])$ satisfy the quantum Yang-Baxter
equation (\ref{qybe}) is that
\begin{eqnarray}
&&r^{1,2}[h]r^{1,3}[h]+r^{1,2}[h]r^{2,3}[h]+r^{1,3}[h]r^{2,3}[h]
+r^{1,2}[h]r^{1,3}[h]r^{2,3}[h]
\nonumber \\
&=&r^{1,3}[h]r^{1,2}[h]+r^{2,3}[h]r^{1,2}[h]+r^{2,3}[h]r^{1,3}[h]
+r^{2,3}[h]r^{1,3}[h]r^{1,2}[h].
\label{condition}
\end{eqnarray}
in the associative algebra $\left( \mathcal{L}\otimes \mathcal{L}\otimes
\mathcal{L}\right) [[h]]$.
\end{theorem}

If $\mathcal{L}^{\prime }$ is the algebra got by appending a unit element to
$\mathcal{L}$ we see by adding the appropriate unit element to both sides
that (\ref{condition}) is equivalent to the condition that $1+r[h]$
satisfies the quantum Yang Baxter equation in $\left( \mathcal{L}^{\prime
}\otimes \mathcal{L}^{\prime }\otimes \mathcal{L}^{\prime }\right) [[h]]$,
\begin{equation}
\left( 1+r^{1,2}[h]\right) \left( 1+r^{1,3}[h]\right) \left(
1+r^{2,3}[h]\right) =\left( 1+r^{2,3}[h]\right) \left( 1+r^{1,3}[h]\right)
\left( 1+r^{1,2}[h]\right) .  \label{toy}
\end{equation}

In \cite{HuPu1} it was claimed erroneously that the \textit{symmetrised}
double product integral $\prod \prod (1+dr[h])$ satisfied the quantum
Yang-Baxter equation if $r[h]$ satisfies (\ref{condition}); in fact for
symmetrised integrals (\ref{condition}) is sufficient for the quantum
Yang-Baxter equation to hold only to order $h^3.$

The condition (\ref{condition}) also occurs in a somewhat different context
in \cite{Enri}. Equating coefficients of $h^2,$ $h^3,\ldots $ in (\ref
{condition}) we find firstly that the first order coefficient $r_1\in
\mathcal{L}\otimes \mathcal{L}$ of $r$ satisfies the classical Yang-Baxter
equation
\begin{equation}
\lbrack r_1^{1,2},r_1^{1,3}]+[r_1^{1,2},r_1^{2,3}]+[r_1^{1,3},r_1^{2,3}]=0,
\label{cybe}
\end{equation}
while for higher order coefficents
\begin{eqnarray}
&&\sum_{s+t=N+1}\left(
[r_s^{1,2},r_t^{1,3}]+[r_s^{1,2},r_t^{2,3}]+[r_s^{1,3},r_t^{2,3}]\right)
\nonumber \\
+ &&\sum_{s+t+u=N+1}\left(
r_s^{1,2}r_t^{1,3}r\,_u^{2,3}-r_u^{2,3}r_t^{1,3}r_s^{1,2}\right) =0
\label{latter}
\end{eqnarray}
where the sums are over ordered pairs $(s,t)$ and triples $(s,t,u)$ of
natural numbers whose sums are $N+1.$ Assuming that $r_1,r_2,\ldots ,r_{N-1}$
have been found, (\ref{latter}) gives a finite dimensional inhomogeneous
linear equation for $r_N.$ It is shown in \cite{Enri} that, for a given
solution $r_1$ of (\ref{cybe}), the resulting hierarchy of equations has
solutions, however, to make this inference, the existence of a quantisation
of an arbitrary Lie bialgebra, as proved in \cite{EtKa}, is needed.

Our second Theorem is that the directed double product integral $\stackrel{%
\rightarrow \leftarrow }{\prod }(1+dr[h])$ has a multiplicative inverse
which is also a directed product integral but with arrows reversed.

\begin{theorem}
For arbitrary $r[h]\in h\left( \mathcal{L}\otimes \mathcal{L}\right) [[h]],$
$\stackrel{\rightarrow \leftarrow }{\prod }(1+dr[h])$ is invertible in $%
\left( \mathcal{T(L)}\otimes \mathcal{T(L)}\right) [[h]]$ with inverse $%
\stackrel{\leftarrow \rightarrow }{\prod }(1+dr^{\prime }[h])$ where $%
r^{\prime }[h]$ is the quasiinverse of $r[h]$ in $h\mathcal{L}[[h]],$ that
is, the unique element satisfying
\[
r[h]+r^{\prime }[h]+r[h]r^{\prime }[h]=r^{\prime }[h]+r[h]+r^{\prime
}[h]r[h]=0.
\]
\end{theorem}

Thus $1+r^{\prime }[h]$ is the inverse of $1+r[h]$ in $\left( \mathcal{L}%
^{\prime }\otimes \mathcal{L}^{\prime }\right) [[h]].$

Let us now demonstrate how Theorems 1 and 2 permit a simple construction of
a deformed Hopf algebra, which as an associative algebra is a subalgebra of
the It\^{o} shuffle algebra $\mathcal{T(L)}[[h]],$ whose semi-classical
limit is an arbitrary quasitriangular Lie bialgebra with commutator Lie
bracket. This may be compared with the general construction of Enriquez \cite
{Enri}.

We first recall \cite{Huds} that the subspace $\mathcal{S(L)}$ of $\mathcal{%
T(L)}$ consisting of symmetric tensors is isomorphic \textit{as a Hopf
algebra} to the universal enveloping algebra $\mathcal{U(L)}$ of the Lie
algebra $\mathcal{L}$. Let us now consider the algebra $\mathcal{T(L)}%
_0[[h]] $ consisting of formal power series with coefficients in $\mathcal{%
T(L)}$ whose zero-order coefficients are in $\mathcal{S(L)}$. We define a
deformed coproduct $\Delta [h]$ on $\mathcal{T(L)}_0[[h]]$ by
\begin{equation}
\Delta [h](\alpha )=R[h]\Delta (\alpha )R[h]^{-1}  \label{10000}
\end{equation}
where $R[h]=$ $\stackrel{\rightarrow \leftarrow }{\prod }(1+dr[h]),$ and
consequently also its inverse $R[h]^{-1}=\stackrel{\leftarrow \rightarrow }{%
\prod }(1+dr^{\prime }[h]),$ satisfy the quantum Yang-Baxter equation (\ref
{qybe}). Using (\ref{qybe}) together with the quasitriangular relations (\ref
{triangular}) and the reversed quasitriangular relations
\[
\left( \Delta \otimes \func{id}\,_{\mathcal{T(L)}}\right)
R[h]^{-1}=R^{2,3}[h]^{-1}R^{1,3}[h]^{-1},\text{ }\left( \func{id}\,_{%
\mathcal{T(L)}}\otimes \Delta \right)
R[h]^{-1}=R^{1,2}[h]^{-1}R^{1,3}[h]^{-1}
\]
satisfied by $R[h]^{-1}$ and the multiplicativity and coassociativity of the
undeformed coproduct $\Delta ,$ we have
\begin{eqnarray*}
&&\left( \Delta [h]\otimes \func{id}\,_{\mathcal{T(L)}}\right) \Delta
[h](\alpha ) \\
&=&R^{1,2}[h]R^{1,3}[h]R^{2,3}[h]\left( \Delta \otimes \func{id}\,_{\mathcal{%
T(L)}}\right) \Delta (\alpha )R^{2,3}[h]^{-1}R^{1,3}[h]^{-1}R^{1,2}[h]^{-1}
\\
&=&R^{2,3}[h]R^{1,3}[h]R^{1,2}[h]\left( \func{id}\,_{\mathcal{T(L)}}\otimes
\Delta \right) \Delta (\alpha )R^{1,2}[h]^{-1}R^{1,3}[h]^{-1}R^{2,3}[h]^{-1}
\\
&=&\left( \func{id}\,_{\mathcal{T(L)}}\otimes \Delta [h]\right) \Delta
[h](\alpha )
\end{eqnarray*}
so that $\Delta [h]$ is coassociative.

The Lie bialgebra cobracket of the semiclassical limit of this deformation
is given \cite{EtSc} by
\[
\delta (L)=h^{-1}\left( \Delta [h](L)-\Delta [h]^{\func{opp}}(L)\right)
+o(h),\text{ }L\in \mathcal{L}
\]
Using the facts that
\[
\stackrel{\rightarrow \leftarrow }{\prod }(1+dr[h])=1\,_{\mathcal{T(L)}%
}+r[h]+o(h^2),\text{ }\stackrel{\leftarrow \rightarrow }{\prod }%
(1+dr^{\prime }[h])=1\,_{\mathcal{T(L)}}+r^{\prime }[h]+o(h^2)
\]
and $\left( r^{\prime }\right) _1=-r_1$ we find that
\[
\delta (L)=[r_1-\tau _{(2,1)}r_1,L^1+L^2].
\]
Since $r_1$ satisfies the classical Yang-Baxter equation, the Lie bialgebra $%
(\mathcal{L}$,$[.,.],\delta )$ is quasitriangular, and every commutator
quasitriangular Lie bialgebra can be quantised in this way.

In Section 2 we review the theory of ordered double product integrals. Some
properties of the coproduct are reviewed in Section 3 and used to derive
another necessary and sufficient condition that a product integral satisfy
the quantum Yang-Baxter equation. In Section 4 we show that this second
condition reduces to (\ref{condition}) and thus complete the proof of
Theorem 1. In Section 5 we prove Theorem 2. In section 6 we describe an
example of our quantisation procedure in detail.

\textbf{Acknowledgement}. The authors wish to thank the referees, and
especially Benjamin Enriquez, for suggestions and criticisms which have
greatly improved the first draft of this paper.

\section{Ordered double product integrals.}

\medskip For $r[h]\in h\left( \mathcal{L\otimes L}\right) [[h]]$ the ordered
double products $\stackrel{\rightarrow \leftarrow }{\prod }(1+dr[h]),$ $%
\stackrel{\leftarrow \rightarrow }{\prod }(1+dr^{\prime }[h])$ are defined
\cite{HuPu} by
\begin{eqnarray*}
\stackrel{\rightarrow \leftarrow }{\prod }(1+dr[h]) &=&_{\mathcal{T(L)}%
}\prod^{\rightarrow }(1+d(\widehat{\prod^{\leftarrow }}_{\mathcal{L}%
}(1+dr[h]))=\prod^{\leftarrow }\,_{\mathcal{T(L)}}(1+d(_{\mathcal{L}}%
\widehat{\prod^{\rightarrow }}(1+dr[h])), \\
\stackrel{\leftarrow \rightarrow }{\prod }(1+dr^{\prime }[h]) &=&_{\mathcal{%
T(L)}}\prod^{\leftarrow }(1+d(\widehat{\prod^{\rightarrow }}_{\mathcal{L}%
}(1+dr[h]))=\prod^{\rightarrow }\,_{\mathcal{T(L)}}(1+d(_{\mathcal{L}}%
\widehat{\prod^{\leftarrow }}(1+dr[h])).
\end{eqnarray*}
Here the notations are as follows. Given a unital associative algebra $%
\mathcal{A}$ and an element $l[h]$ of $h\left( \mathcal{A}\otimes \mathcal{L}%
\right) [[h]]$ the product integrals $_{\mathcal{A}}\stackrel{\rightarrow }{%
\prod }(1+dl[h])$ and $_{\mathcal{A}}\stackrel{\leftarrow }{\prod }(1+dl[h])$
are the unique elements $Y[h]$ and $Z[h]$ of $\left( \mathcal{A}\otimes
\mathcal{T(L)}\right) [[h]]$ which satisfy the left and right differential
equations
\[
\left( \func{id}\,_{\mathcal{A}}\otimes \overrightarrow{d}\right)
Y[h]=Y^{12}[h]l^{13}[h],\left( \func{id}\,_{\mathcal{A}}\otimes \;%
\overleftarrow{d}\right) Z[h]=Z^{13}[h]l^{12}[h],
\]
together with the initial conditions
\[
\left( \func{id}\,_{\mathcal{A}}\otimes \varepsilon \right) Y[h]=1\,_{%
\mathcal{A}},\text{ }\left( \func{id}\,_{\mathcal{A}}\otimes \varepsilon
\right) Z[h]=1\,_{\mathcal{A}}
\]
where $\varepsilon $ is the counit of the It\^{o} Hopf algebra $\mathcal{T(L)%
}$. Equivalently
\begin{equation}
Y[h]=1\,_{\mathcal{A\otimes T(L)}}+\sum_{n=1}^\infty \left(
l^{12}[h]l^{13}[h]\cdots l^{1n+1}[h]\right) ,  \label{unital}
\end{equation}
\[
Z[h]=1\,_{\mathcal{A\otimes T(L)}}+\sum_{n=1}^\infty \left(
l^{1n+1}[h]l^{1n}[h]\cdots l^{12}[h]\right) ,
\]
rearranged as formal power series with coefficients in $\mathcal{A\otimes
T(L)}$. Product integrals $\stackrel{\rightarrow }{\prod }\,_{\mathcal{A}%
}(1+dm[h])$ and $\stackrel{\leftarrow }{\prod }\,_{\mathcal{A}}(1+dm[h]),$%
where now $m[h]\in h$ $\left( \mathcal{L}\otimes \mathcal{A}\right) [[h]],$
are defined analogously in $\left( \mathcal{T(L)\otimes A}\right) [[h]].$ If
$\mathcal{A}$ is not necessarily unital, the decapitated product integrals $%
_{\mathcal{A}}\widehat{\stackrel{\rightarrow }{\prod }}(1+dl[h]),$ $_{%
\mathcal{A}}\widehat{\stackrel{\leftarrow }{\prod }}(1+dl[h]),$ $\,\widehat{%
\stackrel{\rightarrow }{\prod }}_{\mathcal{A}}(1+dm[h]),$ $\widehat{%
\stackrel{\leftarrow }{\prod }}\,_{\mathcal{A}}(1+dm[h])$ are defined by
omitting the initial unit term from the expansion (\ref{unital}) and its
analogues. Thus $\widehat{\stackrel{\leftarrow }{\prod }}_{\mathcal{L}%
}(1+dr[h])$ and $_{\mathcal{L}}\widehat{\stackrel{\rightarrow }{\prod }}%
(1+dr[h])$ are defined as elements of $h\left( \mathcal{T(L)}\otimes
\mathcal{L}\right) [[h]]$ and $h\left( \mathcal{L}\otimes \mathcal{T(L)}%
\right) [[h]]$ respectively, so that $_{\mathcal{T(L)}}\stackrel{%
\rightarrow }{\prod }(1+d(\widehat{\stackrel{\leftarrow }{\prod }}_{\mathcal{%
L}}(1+dr[h]))$ and $\stackrel{\leftarrow }{\prod }\,_{\mathcal{T(L)}}(1+d(_{%
\mathcal{L}}\widehat{\stackrel{\rightarrow }{\prod }}(1+dr[h]))$ are defined
as elements of $\left( \mathcal{T(L)}\otimes \mathcal{T(L)}\right) [[h]].$
It is shown in \cite{HuPu} that these two iterated products are in fact
equal. The double product integral $\stackrel{\leftarrow \rightarrow }{\prod
}(1+dr[h])$ is defined analogously by reversing the directions of all arrows.

Below we shall use the fact, proved in \cite{HuPu}, that $\stackrel{%
\rightarrow \leftarrow }{\prod }(1+dr[h])$ and $\stackrel{\leftarrow
\rightarrow }{\prod }(1+dr[h])$ are both of form
\begin{equation}
1\,_{\mathcal{T(L)\otimes T(L)}}+r[h]+\beta [h]  \label{expansion}
\end{equation}
where $\beta [h]=\sum \beta _{m,n}[h]$ in $\left( \mathcal{T(L)\otimes T(L)}%
\right) [[h]]\mathcal{=}\left( \bigoplus_{m,n=0}^\infty \left( \bigotimes^m%
\mathcal{L}\right) \mathcal{\otimes }\left( \bigotimes^n\mathcal{L}\right)
\right) [[h]]$ has nonzero homogeneous components $\beta _{m,n}[h]$ only if
both $m$ and $n$ are nonzero and at least one exceeds $1.$

\section{Properties of $\Delta .$}

\medskip \medskip In \cite{HuPu} it is shown that the element $R[h]=%
\stackrel{\rightarrow \leftarrow }{\prod }(1+dr[h])$ of $\left( \mathcal{%
T(L)\otimes T(L)}\right) [[h]]$ satisfies the quasitriangularity identities (%
\ref{triangular}) and consequently, for arbitrary $m,n=1,2,\ldots $%
\begin{equation}
\left( \Delta ^{(m)}\otimes \Delta ^{(n)}\right) R[h]=\prod_{(j,k)\in \Bbb{N}%
_m\times \Bbb{N}_n}R^{j,m+n+1-k}[h]  \label{hippo}
\end{equation}
in $\left( \left( \bigotimes^m\mathcal{T(L)}\right) \otimes \left(
\bigotimes^n\mathcal{T(L)}\right) \right) [[h]].$ Here the notation is as
follows. For $m=1,2,\ldots ,\Bbb{N}_m$ denotes the ordered set $(1,2,\ldots
,m),$ and $\Delta ^{(m)}:\mathcal{T(L)}\rightarrow \bigotimes^m\mathcal{T(L)}
$ is the $m$-fold coproduct which may be defined inductively by
\[
\Delta ^{(1)}=\func{id}\,_{\mathcal{T(L)}},\text{ }\Delta ^{(m)}=\left(
\Delta \otimes \func{id}\,_{\otimes ^{m-1}\mathcal{T(L)}}\right) \circ
\Delta ^{(m-1)}.
\]
In particular $\Delta ^{(2)}=\Delta .$ $\stackrel{\rightarrow \leftarrow }{%
\prod }_{(j,k)\in \Bbb{N}_m\times \Bbb{N}_n}R^{j,m+n+1-k}[h]$ is the element
of $\left( \left( \bigotimes^m\mathcal{T(L)}\right) \otimes \left(
\bigotimes^n\mathcal{T(L)}\right) \right) [[h]]$ defined by by the
equivalent prescriptions
\begin{equation}
\prod_{(j,k)\in \Bbb{N}_m\times \Bbb{N}_n}R^{j,m+n+1-k}[h]=\prod_{j=1}^m%
\left\{ \prod_{k=1}^nR^{j,m+n+1-k}[h]\right\} =\prod_{k=1}^n\left\{
\prod_{j=1}^mR^{j,m+n+1-k}[h]\right\} .  \label{prescriptions}
\end{equation}
\medskip

\medskip From the defining property of $\Delta ,$%
\[
\Delta (L_1\otimes L_2\otimes \cdots \otimes L_m)=\sum_{j=0}^m(L_1\otimes
L_2\otimes \cdots \otimes L_j)\otimes (L_{j+1}\otimes L_{j+2}\otimes \cdots
\otimes L_m)
\]
for arbitrary $L_1,L_2,\cdots ,L_m\in \mathcal{L}$, we deduce that for the $%
m $th order coproduct,
\begin{eqnarray*}
&&\Delta ^{(m)}(L_1\otimes L_2\otimes \cdots \otimes L_m) \\
&=&\sum_{0\leq j_1\leq j_2\leq \cdots \leq j_{m-1}\leq m}\left( L_1\otimes
L_2\otimes \cdots \otimes L_{j_1}\right) \otimes \cdots \otimes \left(
L_{j_{m-1}+1}\otimes L_{j_{m-1}+2}\otimes \cdots \otimes L_m\right) \\
&=&L_1\otimes L_2\otimes \cdots \otimes L_m+\beta
\end{eqnarray*}
where $L_1\otimes L_2\otimes \cdots \otimes L_m$ is regarded as an element
of $\bigotimes^m\mathcal{T(L)=}\bigotimes^m\left( \mathcal{\oplus }%
_{n=0}^\infty \left( \bigotimes^n\mathcal{L}\right) \right) $ and the
remainder $\beta $ consists of a sum of elements of each of which is of rank
$>1$ in at least one copy of $\mathcal{\oplus }_{n=0}^\infty \left(
\bigotimes^n\mathcal{L}\right) $. Thus, by linearity, for general $\alpha
=(\alpha _0,\alpha _1,\alpha _2,\ldots )\in \mathcal{T(L)}$,
\begin{equation}
\alpha _m=\left( \Delta ^{(m)}(\alpha )\right) _{(1,1,\ldots ,\stackrel{(m)}{%
1})}  \label{formula}
\end{equation}
in $\bigotimes^m\mathcal{T(L)=}\bigotimes^m\left( \mathcal{\oplus }%
_{n=0}^\infty \left( \bigotimes^n\mathcal{L}\right) \right) .$ (\ref{formula}%
) also holds when $m=0$ if $\Delta ^{(0)}$ is defined to be the counit in $%
\mathcal{T(L)}$.

\medskip Now consider the products $R^{1,2}[h]R^{1,3}[h]R^{2,3}[h]$ and $%
R^{2,3}[h]R^{1,3}[h]R^{1,2}[h]$ in
\[
\left( \mathcal{T(L)\otimes T(L)\otimes T(L)}\right) [[h]]\mathcal{=}\left(
\bigoplus_{m,n,p=0}^\infty (\otimes ^m\mathcal{L)\otimes }(\otimes ^n%
\mathcal{L)\otimes }(\otimes ^p\mathcal{L)}\right) [[h]]
\]
and denote by $\gamma _{m,n,p}[h],$ $\tilde{\gamma}_{m,n,p}[h]$ their
components of joint rank $(m,n,p).$ If any of the nonnegative integers $%
m,n,p $ is zero then $\gamma _{m,n,p}[h]$ and $\tilde{\gamma}_{m,n,p}[h]$
are formed entirely from only one of the three terms $R^{1,2}[h]$, $%
R^{1,3}[h]$ and $R^{2,3}[h]$ and so must be equal. Assume that $(m,n,p)\in
\Bbb{N}^3.$ Then $\gamma _{m,n,p}[h]$ can be characterised as the component
of $\left( \Delta ^{(m)}\otimes \Delta ^{(n)}\otimes \Delta ^{(p)}\right)
\left( R^{12}[h]R^{13}[h]R^{23}[h]\right) $ in $\left( \mathcal{T(L)\otimes
T(L)\otimes T(L)}\right) [[h]]$ of joint rank $((1,1,\ldots ,\stackrel{(m)}{1%
}),(1,1,\ldots ,\stackrel{(n)}{1}),(1,1,\ldots ,\stackrel{(p)}{1})).$ By (%
\ref{hippo}), and the multiplicativity of the iterated coproducts,
\begin{eqnarray*}
&&\left( \Delta ^{(m)}\otimes \Delta ^{(n)}\otimes \Delta ^{(p)}\right)
\left( R^{1,2}[h]R^{1,3}[h]R^{2,3}[h]\right) \\
&=&\prod_{(j,k)\in \Bbb{N}_m\times \Bbb{N}_n}R^{j,m+n+1-k}[h]\prod_{(j,l)\in
\Bbb{N}_m\times \Bbb{N}_p}R^{j,m+n+p+1-l}[h]\prod_{(k,l)\in \Bbb{N}_n\times
\Bbb{N}_p}R^{m+k,m+n+p+1-l}[h].
\end{eqnarray*}
Thus $\gamma _{m,n,p}[h]$ is the component of the same joint rank in this
expression. By (\ref{expansion}) this is the component of the same joint
rank in
\begin{eqnarray}
&&\prod_{(j,k)\in \Bbb{N}_m\times \Bbb{N}_n}\left( 1\,_{\mathcal{T(L)\otimes
T(L)}}+r[h]^{j,m+n+1-k}\right) \prod_{(j,l)\in \Bbb{N}_m\times \Bbb{N}%
_p}\left( 1\,_{\mathcal{T(L)\otimes T(L)}}+r[h]^{j,m+n+p+1-l}\right)
\nonumber \\
&&\prod_{(k,l)\in \Bbb{N}_n\times \Bbb{N}_p}\left( 1\,_{\mathcal{T(L)\otimes
T(L)}}+r[h]^{m+k,m+n+p+1-l}\right)  \label{p1}
\end{eqnarray}
since higher order terms than $r[h]$ in the expansion (\ref{expansion}) can
only contribute to components in $\mathcal{T(L)\otimes T(L)=}%
\bigoplus_{m,n=0}^\infty (\otimes ^m\mathcal{L)\otimes }(\otimes ^n\mathcal{%
L)}$ at least one of whose ranks exceeds $1.$ By expressing $\tilde{\gamma}%
_{m,n,p}[h]$ similarly as the component of the same joint rank in
\begin{eqnarray}
&&\prod_{(k,l)\in \Bbb{N}_n\times \Bbb{N}_p}\left( 1\,_{\mathcal{T(L)\otimes
T(L)}}+r[h]^{m+k,m+n+p+1-k}\right) \prod_{(j,l)\in \Bbb{N}_m\times \Bbb{N}%
_p}\left( 1\,_{\mathcal{T(L)\otimes T(L)}}+r[h]^{j,m+n+p+1-l}\right)
\nonumber \\
&&\prod_{(k,l)\in \Bbb{N}_n\times \Bbb{N}_p}\left( 1\,_{\mathcal{T(L)\otimes
T(L)}}+r[h]^{j,m+n+1-l}\right)  \label{p2}
\end{eqnarray}
we complete the proof of the following.

\begin{theorem}
A necessary and sufficient condition that $R[h]=$ $\stackrel{\rightarrow
\leftarrow }{\prod }(1+dr[h])$ satisfy the quantum Yang-Baxter equation is
that, for all triples $(m,n,p)$ of natural numbers, the components of joint
rank $((1,1,\ldots ,\stackrel{(m)}{1}),(1,1,\ldots ,\stackrel{(n)}{1}%
),(1,1,\ldots ,\stackrel{(p)}{1}))$ in the products (\ref{p1}) and (\ref{p2}%
) are equal.
\end{theorem}

\section{Proof of Theorem 1.}

From Theorem 3 we deduce that a necessary and sufficient condition for $R[h]=%
\stackrel{\rightarrow \leftarrow }{\prod }(1+dr[h])$ to satisfy the quantum
Yang-Baxter equation is that, for all $(m,n,p)\in \Bbb{N}^3,$%
\begin{eqnarray}
&&\prod_{(j,k)\in \Bbb{N}_m\times \Bbb{N}_n}\rho ^{j,m+n+1-k}\prod_{(j,l)\in
\Bbb{N}_m\times \Bbb{N}_p}\rho ^{j,m+n+p+1-l}\prod_{(k,l)\in \Bbb{N}_n\times
\Bbb{N}_p}\rho ^{m+k,m+n+p+1-l}  \nonumber \\
&=&\prod_{(k,l)\in \Bbb{N}_n\times \Bbb{N}_p}\rho
^{m+k,m+n+p+1-l}\prod_{(j,l)\in \Bbb{N}_m\times \Bbb{N}_p}\rho
^{j,m+n+p+1-l}\prod_{(j,k)\in \Bbb{N}_m\times \Bbb{N}_n}\rho ^{j,m+n+1-k}
\label{braces1}
\end{eqnarray}
in the associative algebra $\left( \left( \bigotimes^m\mathcal{L}^{\prime
}\right) \mathcal{\otimes }\left( \bigotimes^n\mathcal{L}^{\prime }\right)
\mathcal{\otimes }\left( \bigotimes^p\mathcal{L}^{\prime }\right) \right)
[[h]],$ where $\rho =1+r[h]$ and $h$ is suppressed for brevity. Indeed it is
precisely these ``pure It\^{o}'' terms which can contribute to the
components of joint rank $((1,1,\ldots ,\stackrel{(m)}{1}),(1,1,\ldots ,%
\stackrel{(n)}{1}),(1,1,\ldots ,\stackrel{(p)}{1}))$ of the multitensors (%
\ref{p1}) and (\ref{p2}).

Taking $m=n=p=1$ we find that it is necessary that
\begin{equation}
\rho ^{1,2}\rho ^{1,3}\rho ^{2,3}=\rho ^{2,3}\rho ^{1,3}\rho ^{1,2}
\end{equation}
from which we deduce (\ref{condition}) by expanding and cancelling the unit
on each side.\medskip

\medskip To prove the sufficiency of (\ref{condition}), or equivalently of
(16), we must show that it implies (\ref{braces1}) for arbitrary $(m,n,p)\in
\Bbb{N}^3.$ \medskip First we prove the case when $m=n=1$ by induction on $%
p; $ when $p=1$ this is just (16). With $m=n=1$ the left hand side of (\ref
{braces1}) becomes
\begin{eqnarray*}
&\medskip &\rho ^{1,2}\left\{ \rho ^{1,p+2}\rho ^{1,p+1}\cdots \rho
^{1,3}\right\} \left\{ \rho ^{2,p+2}\rho ^{2,p+1}\cdots \rho ^{2,3}\right\}
\\
&=&\left( \rho ^{1,2}\rho ^{1,p+2}\rho ^{2,p+2}\right) \left\{ \rho
^{1,p+1}\rho ^{1,p}\cdots \rho ^{1,3}\right\} \left\{ \rho ^{2,p+1}\rho
^{2,p}\cdots \rho ^{2,3}\right\} \\
&=&\rho ^{2,p+2}\rho ^{1,p+2}\rho ^{1,2}\left\{ \rho ^{1,p+1}\rho
^{1,p}\cdots \rho ^{1,3}\right\} \left\{ \rho ^{2,p+1}\rho ^{2,p}\cdots \rho
^{2,3}\right\}
\end{eqnarray*}
using commutativity of $\rho ^{2,p+2}$ with all but the first element of $%
\left\{ \rho ^{1,p+2}\rho ^{1,p+1}\cdots \rho ^{1,3}\right\} $ and (16). By
the inductive hypothesis this is
\begin{eqnarray*}
&&\rho ^{2,p+2}\rho ^{1,p+2}\left\{ \rho ^{2,p+1}\rho ^{2,p}\cdots \rho
^{2,3}\right\} \left\{ \rho ^{1,p+1}\rho ^{1,p}\cdots \rho ^{1,3}\right\}
\rho ^{1,2} \\
&=&\left\{ \rho ^{2,p+2}\rho ^{2,p+1}\cdots \rho ^{2,3}\right\} \left\{ \rho
^{1,p+2}\rho ^{1,p+1}\cdots \rho ^{1,3}\right\} \rho ^{1,2}
\end{eqnarray*}
which is the right hand side.

$\medskip $Next we prove the case $n=1$ by induction on $m;$ when $m=1$ this
becomes the case just proved. We have to show that, for arbitrary $m\,$and $%
p,$%
\begin{eqnarray}
&&\prod_{j=1}^m\rho ^{j,m+1}\prod_{(j,l)\in \Bbb{N}_m\times \Bbb{N}_p}\rho
^{j,m+2+p-l}\prod_{l=1}^p\rho ^{m+1,m+p+2-l}  \nonumber \\
&=&\prod_{l=1}^p\rho ^{m+1,m+p+2-l}\prod_{(j,l)\in \Bbb{N}_m\times \Bbb{N}%
_p}\rho ^{j,m+2+p-l}\prod_{j=1}^m\rho ^{j,m+1}  \label{7}
\end{eqnarray}
Using the first definition in (\ref{prescriptions}) to define the double
product and extracting the last internal product to the right, we express
the left hand side as
\begin{eqnarray*}
&&\prod_{j=1}^m\rho ^{j,m+1}\prod_{(j,l)\in \Bbb{N}_{m-1}\times \Bbb{N}%
_p}\rho ^{j,m+2+p-l}\prod_{l=1}^p\rho ^{m,m+p+2-l}\prod_{l=1}^p\rho
^{m+1,m+p+2-l} \\
&=&\prod_{j=1}^{m-1}\rho ^{j,m+1}\prod_{(j,l)\in \Bbb{N}_{m-1}\times \Bbb{N}%
_p}\rho ^{j,m+2+p-l}\left( \rho ^{m,m+1}\rho ^{m,m+p+1}\rho
^{m+1,m+p+1}\right) \\
&&\rho ^{m,m+p}\rho ^{m,m+p-1}\cdots \rho ^{m,m+2}\rho ^{m+1,m+p}\rho
^{m+1,m+p-1}\cdots \rho ^{m+1,m+2}
\end{eqnarray*}
$\medskip $using the commutativity of $\rho ^{m,m+1}$ with all remaining
elements of the double product, and of $\rho ^{m+1,m+p+1}$ with all but the
first element of the extracted internal product. Reversing the bracketed
Yang-Baxter triple, the terms after $\prod_{j=1}^{m-1}\rho
^{j,m+1}\prod_{(j,l)\in \Bbb{N}_{m-1}\times \Bbb{N}_p}\rho ^{j,m+2+p-l}$
become
\begin{eqnarray*}
&&\left( \rho ^{m+1,m+p+1}\rho ^{m,m+p+1}\rho ^{m,m+1}\right) \rho
^{m,m+p}\rho ^{m,m+p-1}\cdots \rho ^{m,m+2} \\
&&\rho ^{m+1,m+p}\rho ^{m+1,m+p-1}\cdots \rho ^{m+1,m+2} \\
&=&\rho ^{m+1,m+p+1}\rho ^{m,m+p+1}\left( \rho ^{m,m+1}\rho ^{m,m+p}\rho
^{m+1,m+p}\right) \rho ^{m,m+p-1}\rho ^{m,m+p-2}\cdots \rho ^{m,m+2} \\
&&\rho ^{m+1,m+p-1}\rho ^{m+1,m+p-2}\cdots \rho ^{m+1,m+2} \\
&=&\rho ^{m+1,m+p+1}\rho ^{m,m+p+1}\left( \rho ^{m+1,m+p}\rho ^{m,m+p}\rho
^{m,m+1}\right) \rho ^{m,m+p-1}\rho ^{m,m+p-2}\cdots \rho ^{m,m+2} \\
&&\rho ^{m+1,m+p-1}\rho ^{m+1,m+p-2}\cdots \rho ^{m+1,m+2} \\
&=&\rho ^{m+1,m+p+1}\rho ^{m,m+p+1}\rho ^{m+1,m+p}\rho ^{m,m+p}\left( \rho
^{m,m+1}\rho ^{m,m+p-1}\rho ^{m+1,m+p-1}\right) \\
&&\rho ^{m,m+p-2}\cdots \rho ^{m,m+2}\rho ^{m+1,m+p-2}\rho
^{m+1,m+p-3}\cdots \rho ^{m+1,m+2} \\
&=&\rho ^{m+1,m+p+1}\rho ^{m,m+p+1}\rho ^{m+1,m+p}\rho ^{m,m+p}\left( \rho
^{m+1,m+p-1}\rho ^{m,m+p-1}\rho ^{m,m+1}\right) \\
&&\rho ^{m,m+p-2}\cdots \rho ^{m,m+2}\rho ^{m+1,m+p-2}\rho
^{m+1,m+p-3}\cdots \rho ^{m+1,m+2} \\
&=&\cdots \\
&=&\rho ^{m+1,m+p+1}\rho ^{m,m+p+1}\rho ^{m+1,m+p}\rho ^{m,m+p}\rho
^{m+1,m+p-1}\rho ^{m,m+p-1}\cdots \\
&&\rho ^{m+1,m+2}\rho ^{m,m+2}\rho ^{m,m+1} \\
&=&\rho ^{m+1,m+p+1}\rho ^{m+1,m+p}\rho ^{m+1,m+p-1}\cdots \rho
^{m+1,m+2}\rho ^{m,m+p+1}\rho ^{m,m+p}\rho ^{m,m+p-1}\cdots \\
&&\rho ^{m,m+2}\rho ^{m,m+1}.
\end{eqnarray*}
Restoring the omitted initial terms $\prod_{j=1}^{m-1}\rho
^{j,m+1}\prod_{(j,l)\in \Bbb{N}_{m-1}\times \Bbb{N}_p}\rho ^{j,m+2+p-l}$ we
obtain
\[
\prod_{j=1}^{m-1}\rho ^{j,m+1}\left( \prod_{(j,l)\in \Bbb{N}_{m-1}\times
\Bbb{N}_p}\rho ^{j,m+2+p-l}\right) \prod_{l=1}^p\rho
^{m+1,m+p+2-l}\prod_{l=1}^p\rho ^{m,m+p+2-l}\rho ^{m,m+1}.
\]
Applying the inductive hypothesis to the first three products, this becomes
\[
\prod_{l=1}^p\rho ^{m+1,m+p+2-l}\left( \prod_{(j,l)\in \Bbb{N}_{m-1}\times
\Bbb{N}_p}\rho ^{j,m+2+p-l}\right) \prod_{j=1}^{m-1}\rho
^{j,m+1}\prod_{l=1}^p\rho ^{m,m+p+2-l}\rho ^{m,m+1}
\]
\begin{eqnarray*}
&=&\prod_{l=1}^p\rho ^{m+1,m+p+2-l}\left( \prod_{(j,l)\in \Bbb{N}%
_{m-1}\times \Bbb{N}_p}\rho ^{j,m+2+p-l}\right) \prod_{l=1}^p\rho
^{m,m+p+2-l}\prod_{j=1}^m\rho ^{j,m+1} \\
&=&\prod_{l=1}^p\rho ^{m+1,m+p+2-l}\left( \prod_{(j,l)\in \Bbb{N}_m\times
\Bbb{N}_p}\rho ^{j,m+2+p-l}\right) \prod_{j=1}^m\rho ^{j,m+1}
\end{eqnarray*}
as required, where the last internal product of the original iterated
product is restored again using the first definition in (\ref{prescriptions}%
) to define the corresponding double product.

\medskip Finally we prove the general case by induction on $n$ ; when $n=1$
this is the case just proved.

We shall evaluate the first and last double products on the left hand side
of (\ref{braces1}) using respectively the second and first definitions of (%
\ref{prescriptions}). Detaching the last internal product from the first
double product to the right and the first internal product from the third
double product to the left, we obtain
\begin{eqnarray*}
&&\left( \prod_{k=1}^{n-1}\left\{ \prod_{j=1}^m\rho ^{j,m+n+1-k}\right\}
\right) \prod_{j=1}^m\rho ^{j,m+1}\prod_{(j,l)\in \Bbb{N}_m\times \Bbb{N}%
_p}\rho ^{j,m+n+p+1-l}\prod_{l=1}^p\rho ^{m+1,m+n+p+1-l} \\
&&\left( \prod_{k=2}^n\left\{ \prod_{l=1}^p\rho ^{m+k,m+n+p+1-l}\right\}
\right)
\end{eqnarray*}
The middle three products in this expression can be reversed by (\ref{7}) to
get
\begin{eqnarray}
&&\left( \prod_{k=1}^{n-1}\left\{ \prod_{j=1}^m\rho ^{j,m+n+1-k}\right\}
\right) \prod_{l=1}^p\rho ^{m+1,m+n+p+1-l}\prod_{(j,l)\in \Bbb{N}_m\times
\Bbb{N}_p}\rho ^{j,m+n+p+1-l}\prod_{j=1}^m\rho ^{j,m+1}  \nonumber \\
&&\left( \prod_{k=2}^n\left\{ \prod_{l=1}^p\rho ^{m+k,m+n+p+1-l}\right\}
\right)  \label{grand}
\end{eqnarray}
Since the elements $\rho ^{m+1,m+n+p},\rho ^{m+1,m+n+p-1},\cdots ,\rho
^{m+1,m+n+1}$ commute with those of the first iterated product in (\ref
{grand}) and the elements $\rho ^{1,m+1},\rho ^{2,m+1},\cdots ,\rho ^{m,m+1}$
commute with those of the final iterated product, this is equal to
\begin{eqnarray*}
&&\prod_{l=1}^p\rho ^{m+1,m+n+p+1-l}\left( \prod_{k=1}^{n-1}\left\{
\prod_{j=1}^m\rho ^{j,m+n+1-k}\right\} \right) \prod_{(j,l)\in \Bbb{N}%
_m\times \Bbb{N}_p}\rho ^{j,m+n+p+1-l} \\
&&\left( \prod_{k=2}^n\left\{ \prod_{l=1}^p\rho ^{m+k,m+n+p+1-l}\right\}
\right) \prod_{j=1}^m\rho ^{j,m+1}
\end{eqnarray*}
We now use the inductive hypothesis to reverse the middle three double
products to get
\begin{eqnarray*}
&&\prod_{l=1}^p\rho ^{m+1,m+n+p+1-l}\left( \prod_{k=2}^n\left\{
\prod_{l=1}^p\rho ^{m+k,m+n+p+1-l}\right\} \right) \prod_{(j,l)\in \Bbb{N}%
_m\times \Bbb{N}_p}\rho ^{j,m+n+p+1-l} \\
&&\left( \prod_{k=1}^{n-1}\left\{ \prod_{j=1}^m\rho ^{j,m+n+1-k}\right\}
\right) \prod_{j=1}^m\rho ^{j,m+1}
\end{eqnarray*}
Finally we reabsorb the first and last products into the iterated products,
to obtain
\begin{eqnarray*}
&&\left( \prod_{k=1}^n\left\{ \prod_{l=1}^p\rho ^{m+k,m+n+p+1-l}\right\}
\right) \prod_{(j,l)\in \Bbb{N}_m\times \Bbb{N}_p}\rho ^{j,m+n+p+1-l}\left(
\prod_{k=1}^n\left\{ \prod_{j=1}^m\rho ^{j,m+n+1-k}\right\} \right) \\
&=&\prod_{(k,l)\in \Bbb{N}_n\times \Bbb{N}_p}\rho
^{m+k,m+n+p+1-l}\prod_{(j,l)\in \Bbb{N}_m\times \Bbb{N}_p}\rho
^{j,m+n+p+1-l}\prod_{(j,k)\in \Bbb{N}_m\times \Bbb{N}_n}\rho ^{j,m+n+1-k}
\end{eqnarray*}
as required. $\square $

\section{Proof of Theorem 2.}

Being of the form $1\,_{\mathcal{T(L)\otimes T(L)}}+o(h),$ $R[h]=\stackrel{%
\rightarrow \leftarrow }{\prod }(1+dr[h])$ possesses a multiplicative
inverse $R[h]^{-1}$ of the same form. By the multiplicativity of $\Delta $
and the quasitriangularity relation $\left( \Delta \otimes \func{id}\,_{%
\mathcal{T(L)}}\right) R[h]=R^{13}[h]R^{23}[h]$
\begin{eqnarray*}
1\,_{\mathcal{T(L)\otimes T(L)\otimes T(L)}} &=&\left( \Delta \otimes \func{%
id}\,_{\mathcal{T(L)}}\right) 1\,_{\mathcal{T(L)\otimes T(L)}}=\left( \Delta
\otimes \func{id}\,_{\mathcal{T(L)}}\right) \left( R[h]R[h]^{-1}\right) \\
&=&R^{1,3}[h]R^{2,3}[h]\left( \Delta \otimes \func{id}\,_{\mathcal{T(L)}%
}\right) \left( R[h]^{-1}\right) .
\end{eqnarray*}
On the other hand the inverse of the product $R^{1,3}[h]R^{2,3}[h]$ is the
product of the inverses in the reverse order $\left( R^{2,3}[h]\right)
^{-1}\left( R^{1,3}[h]\right) ^{-1}.$ By uniqueness of inverses it follows
that
\[
\left( \Delta \otimes \func{id}\,_{\mathcal{T(L)}}\right) \left(
R[h]^{-1}\right) =\left( R^{2,3}[h]\right) ^{-1}\left( R^{1,3}[h]\right)
^{-1}.
\]
A similar argument shows that
\[
\left( \func{id}\,_{\mathcal{T(L)}}\otimes \Delta \right) \left(
R[h]^{-1}\right) =\left( R^{1,2}[h]\right) ^{-1}\left( R^{1,3}[h]\right)
^{-1}.
\]
Hence $R[h]^{-1}$ satisfies the reverse quasitriangularity relations. Since
it is manifestly nonzero, by the characterisation theorem of \cite{HuPu} it
must be of form $R[h]^{-1}=\stackrel{\leftarrow \rightarrow }{\prod }%
(1+dr^{\prime }[h])$ for some $r^{\prime }[h]\in h\left( \mathcal{L\otimes L}%
\right) [[h]].$ By comparing terms of rank $\leq 1$ in each copy of $%
\mathcal{T(L)}$ in the identity $R[h]R[h]^{-1}=1$ we find that $r^{\prime
}[h]$ is the quasiinverse of $r[h].$\medskip $\square $

\section{An example.}

We consider the two-dimensional complex associative algebra $\mathcal{L}$
spanned by elements $L,K$ with multiplication $LK=L,$ $K^2=K,$ $L^2=KL=0$
which has the matrix representation
\[
L\mapsto \left[
\begin{array}{cc}
0 & 1 \\
0 & 0
\end{array}
\right] ,\text{ }K\mapsto \left[
\begin{array}{cc}
0 & 0 \\
0 & 1
\end{array}
\right] ,
\]
and for which the commutator Lie algebra is the non-Abelian two-dimensional
Lie algebra. Then $\mathcal{L\wedge L}$ is spanned by the element $%
r_1=L\otimes K-K\otimes L.$ $r_1$ satisfies the classical Yang-Baxter
equation; it defines the unique coboundary Lie bialgebra structure \cite
{EtSc} on this Lie algebra. Also $r_1^2=0$ in the associative algebra $%
\mathcal{L\otimes L}$ and
\begin{equation}
r_1^{1,2}r_1^{1,3}r_1^{2,3}=r_1^{2,3}r_1^{1,3}r_1^{1,2}=0  \label{last}
\end{equation}
in $\mathcal{L\otimes L\otimes L}$. It follows from (\ref{last}) that the
equations (\ref{latter}) can be satisfied by taking all $r_N$ with $N>1$ to
be $0.$ Thus by Theorem 1, $R[h]=\stackrel{\rightarrow \leftarrow }{\prod }%
(1+hdr_1)$ satisfies the quantum Yang-Baxter equation. Since $r_1^2=0,$ the
quasiinverse of $hr_1$ in $h\mathcal{L}[[h]]$ is $-hr_1$ so that by Theorem
2 the inverse of $R[h]$ in $\left( \mathcal{T(L)\otimes T(L)}\right) [[h]]$
is $R[h]^{-1}=\stackrel{\leftarrow \rightarrow }{\prod }(1-hdr_1).$

\medskip To find the actions $\Delta [h](L)$ and $\Delta [h](K)$ of the
deformed coproduct $\Delta [h]$ on $L$ and $K,$ we use the
commutation
relations in $\mathcal{L\otimes L}$%
\begin{equation}
\lbrack r_1,L^1]=L\otimes L,\text{ }[r_1,L^2]=-L\otimes L,\text{ }[r_1,K^1%
]=L\otimes K,\text{ }[r_1,K^2]=-K\otimes L.  \label{commutator4}
\end{equation}
Then, by (\ref{formula}), to find the homogeneous component of joint rank $%
(m,n)$ of
\[
\Delta [h](L)=\stackrel{\rightarrow \leftarrow }{\prod }(1+hdr_1)\Delta (L)%
\stackrel{\leftarrow \rightarrow }{\prod }(1-hdr_1)=\stackrel{\rightarrow
\leftarrow }{\prod }(1+hdr_1)(L^1+L^2)\stackrel{\leftarrow \rightarrow }{%
\prod }(1-hdr_1)
\]
in $\left( \mathcal{T(L)\otimes T(L)}\right) [[h]]$, we compute the
homogeneous component of joint rank $(1,1,...,\stackrel{(m)}{1},1,1,...,%
\stackrel{(n)}{1})$ of
\begin{eqnarray*}
&&\left( \Delta ^{(m)}\otimes \Delta ^{(n)}\right) \Delta [h](L) \\
&=&\left( \Delta ^{(m)}\otimes \Delta ^{(n)}\right) \left( \stackrel{%
\rightarrow \leftarrow }{\prod }(1+hdr_1)\right) \left( \Delta ^{(m)}\otimes
\Delta ^{(n)}\right) (L^1+L^2) \\
&&\left( \Delta ^{(m)}\otimes \Delta ^{(n)}\right) \left( \stackrel{%
\leftarrow \rightarrow }{\prod }(1-hdr_1)\right) \\
&=&\prod_{(j,k)\in \Bbb{N}_m\times \Bbb{N}_n}\left( \stackrel{\rightarrow
\leftarrow }{\prod }(1+hdr_1)\right) ^{j,m+n+1-k}(L^1+L^2+\cdots +L^{m+n}) \\
&&\prod_{(j,k)\in \Bbb{N}_m\times \Bbb{N}_n}\left( \stackrel{\leftarrow
\rightarrow }{\prod }(1-hdr_1)\right) ^{m+1-j,m+k}
\end{eqnarray*}
\medskip in $\left( \otimes ^m\mathcal{T(L)}\right) \mathcal{\otimes }\left(
\otimes ^n\mathcal{T(L)}\right) $, or equivalently of the product in $\left(
\otimes ^m\mathcal{L}^{\prime }\right) \mathcal{\otimes }\left( \otimes ^n%
\mathcal{L}^{\prime }\right) $
\begin{eqnarray*}
&&\prod_{(j,k)\in \Bbb{N}_m\times \Bbb{N}_n}(1+hr_1)^{j,m+n+1-k}(L^1+L^2+%
\cdots +L^{m+n})\prod_{(j,k)\in \Bbb{N}_m\times \Bbb{N}%
_n}(1-hr_1)^{m+1-j,m+k} \\
&=&L^1+L^2+\cdots +L^{m+n}
\end{eqnarray*}
since by starting from the right of the first double product, each term $%
(1+hdr_1)^{j,m+n+1-k}$ in turn commutes with $L^j+L^{m+n+1-k}$ by (\ref
{commutator4}) and with all the other $L^l$ since they occupy disjoint
spaces, and may thus be cancelled with the its inverse in the second double
product. This is evidently zero unless $(m,n)=(1,0)$ or $(0,1).$ It follows
that $\Delta [h](L)=L\otimes 1+1\otimes L.$

Similarly, the $(m,n)$th component of $\Delta [h](K)$ is the $(1,1,...,%
\stackrel{(m)}{1},1,1,...,\stackrel{(n)}{1})$th component of
\[
\prod_{(j,k)\in \Bbb{N}_m\times \Bbb{N}_n}(1+hr_1)^{j,m+n+1-k}(K^1+K^2+%
\cdots +K^{m+n})\prod_{(j,k)\in \Bbb{N}_m\times \Bbb{N}%
_n}(1-hr_1)^{m+1-j,m+k},
\]
whence, using (\ref{commutator4}), we find $(\Delta [h](K))_{0,0}=0,$ $%
(\Delta [h](K))_{1,0}=(\Delta [h](K))_{0,1}=K,$%
\[
(\Delta [h](K))_{1,1}=\left( (1+hr_1)(K^1+K^2)(1-hr_1)\right)
_{1,1}=h(L\otimes K-K\otimes L)
\]
\begin{eqnarray*}
(\Delta [h](K))_{2,1} &=&\left(
(1+hr_1)^{1,3}(1+hr_1)^{2,3}(K^1+K^2+K^3)(1-hr_1)^{2,3}(1-hr_1)^{1,3}\right)
_{1,1,1} \\
&=&\left( (1+hr_1)^{1,3}h(L^2K^3-K^2L^3)(1-hr_1)^{1,3}\right) _{1,1,1} \\
&=&h^2(L\otimes K-K\otimes L)\otimes L
\end{eqnarray*}
and generally for $m\geq 2$%
\[
(\Delta [h](K))_{m,1}=(-h)^m(\otimes ^{m-2}(L))\otimes (L\otimes K-K\otimes
L)\otimes L,
\]
while similarly, for $n\geq 2$%
\[
(\Delta [h](K))_{1,n}=h^nL\otimes (L\otimes K-K\otimes L)\otimes \left(
\otimes ^{n-2}\left( L\right) \right) .
\]
Also, using (\ref{commutator4}) repeatedly and the fact that $hr_1^{1,4}$
commutes with $L^1K^4-K^4L^1=r_1^{1,4}$
\begin{eqnarray*}
(\Delta [h](K))_{2,2} &=&\left(
(1+hr_1)^{1,4}(1+hr_1)^{1,3}(1+hr_1)^{2,4}(1+hr_1)^{2,3}(K^1+K^2+K^3+K^4)%
\right. \\
&&\left. (1-hr_1)^{2,3}(1-hr_1)^{2,4}(1-hr_1)^{1,3}(1-hr_1)^{1,4}\right)
_{1,1,1,1} \\
&=&\cdots =(1+hr_1)^{1,4}h^3\left( -K^1L^2L^3L^4+L^1L^2L^3K^4\right)
(1-hr_1)=0
\end{eqnarray*}
whence also $(\Delta [h](K))_{m,n}=0$ whenever $m,n\geq 2.$ Thus $\Delta
[h](K)$ is completely determined. Note that the components of ranks $(m,1)$
and $(1,n)$ for $m,n\geq 2$ do not belong to $\left( \mathcal{S(L)\otimes
S(L)}\right) [[h]].$

\end{document}